\documentclass[12pt,reqno]{amsart}

\numberwithin{equation}{section}
\newcommand*{\LatexDef}{.}
\usepackage{amsthm}

\usepackage{amsmath}

\usepackage{amssymb}

\usepackage{thmtools}

\usepackage[T1]{fontenc}

\usepackage[shortlabels]{enumitem}
\setlist[enumerate]{leftmargin=*,font=\upshape,align=parleft,label=(\alph*),widest=m}
\setlist[itemize]{leftmargin=*}

\usepackage{xspace}

\usepackage{calrsfs}

\DeclareSymbolFont{defaultmathcal}{OMS}{zplm}{m}{n}
\DeclareSymbolFontAlphabet{\mathcal}{defaultmathcal}
\DeclareSymbolFont{handwritten}{OMS}{rsfs}{m}{n}
\DeclareSymbolFontAlphabet{\handcal}{handwritten}

\usepackage[alphabetic, nobysame]{amsrefs}

\usepackage{hyperref}
\hypersetup{pdfstartview={XYZ null null 1.00}, pdfpagemode=UseNone, colorlinks,breaklinks,linkcolor=blue,urlcolor=blue, anchorcolor=blue,citecolor=blue}

\usepackage[capitalise]{cleveref}

\usepackage{xparse} 

\usepackage{geometry}
\usepackage{microtype}

\newcommand{\N}{\mathbb{N}}

\newcommand{\R}{\mathbb{R}}
\newcommand{\C}{\mathbb{C}}

\newcommand{\Hil}{\mathcal{H}}

\newcommand{\0}{\mathbb{\emptyset}}
\newcommand{\w}{\infty}

\newcommand{\e}{\varepsilon}
\renewcommand{\phi}{\varphi}

\newcommand{\AC}{\mathcal{A}} 
\newcommand{\BC}{\mathcal{B}}
\newcommand{\CC}{\mathcal{C}}
\newcommand{\FC}{\mathcal{F}}

\newcommand{\Pow}{\handcal{P}}
\newcommand{\cl}{\overline}
\newcommand{\Fin}[1][\N]{{#1^{<\N}}}

\newcommand{\actson}{\curvearrowright}

\newcommand*{\defeq}{\mathrel{\vcenter{\baselineskip0.5ex \lineskiplimit0pt \hbox{\scriptsize.}\hbox{\scriptsize.}}}=}

\newcommand*{\defequivlong}{\mathrel{\vcenter{\baselineskip0.5ex \lineskiplimit0pt\hbox{\scriptsize.}\hbox{\scriptsize.}}}\Longleftrightarrow}

\newcommand{\normlr}[1]{\left\|#1\right\|}
\newcommand{\norm}[2][]{#1\| #2 #1\|}

\newcommand{\Los}{\L{}o\'{s}\xspace}
\newcommand{\Slawek}{S\l{}awek\xspace}

\newcommand{\set}[1]{\left\{ #1 \right\}}

\newcommand{\gen}[1]{\left\langle #1 \right\rangle}
\newcommand{\genbig}[2][\big]{#1\langle #2 #1\rangle}

\newcommand{\rest}[1]{\mathord{|_{#1}}}

\newcommand{\eqcomment}[2][6pt]{\Big[\text{#2}\Big] \hspace{#1}}
\newcommand{\eqcommentlr}[2][6pt]{\left[\text{#2}\right] \hspace{#1}}

\newcommand{\onum}[2][th]{$#2^\text{#1}$}

\theoremstyle{plain}

\newtheorem{theorem}[equation]{Theorem}

\crefname{prop}{Proposition}{Propositions}
\newtheorem{prop}[equation]{Proposition}

\newtheorem{lemma}[equation]{Lemma}

\crefname{cor}{Corollary}{Corollaries}
\newtheorem{cor}[equation]{Corollary}

\crefname{obs}{Observation}{Observations}
\newtheorem{obs}[equation]{Observation}

\theoremstyle{definition}

\crefname{defn}{Definition}{Definitions}
\newtheorem{defn}[equation]{Definition}

\crefname{example}{Example}{Examples}
\newtheorem{example}[equation]{Example}

\theoremstyle{remark}

\crefname{remark}{Remark}{Remarks}
\newtheorem{remark}[equation]{Remark}

\newenvironment{remarklike}[2][]{\refstepcounter{equation}\par\medskip\noindent \textit{#2}~\theequation\xspace#1.\rmfamily}{\medskip}
\newenvironment{remarklike*}[2][]{\par\medskip\noindent \textit{#2}#1\textbf{.}\rmfamily\xspace}{\medskip}

\crefname{claim+}{Claim}{Claims}
\newtheorem{claim+}[equation]{Claim}

\newenvironment{case*}[1]{\smallskip\par\noindent \textit{Case}~#1:\rmfamily}{}
\newenvironment{case}[2]{\smallskip\par\noindent \textit{Case}~#1: \rmfamily #2.}{}

\crefname{notation}{Notation}{Notations}
\newtheorem{notation}[equation]{Notation}

\newcommand{\fntsz}[1][11]{\fontsize{#1}{#1}\selectfont}

\newenvironment{acknowledgements}[1][11]{\medskip \fntsz[#1]\begin{trivlist}
		\item[\hskip \labelsep {\textit{Acknowledgements}.}]}{\end{trivlist}\smallskip}

\newenvironment{credits}[1][11]{\smallskip \fntsz[#1]\begin{trivlist}
		\item[\hskip \labelsep {\textit{Credits}.}]}{\end{trivlist}\smallskip}

\newenvironment{enumref}[2][i]
{
	\begin{enumerate}[\upshape(\ref*{#2}.#1)]
}
{
\end{enumerate}
}

\newcommand{\exampleslabelref}[2]{\ref{#1}\ref{#2}\xspace}
\newcommand{\examplesref}[2]{Example \exampleslabelref{#1}{#2}}

\newcommand{\multiexamplesref}[1]{Examples \ref{#1}}

\crefname{examples}{Examples}{Examples}

\newenvironment{examples}[1][\alph*]
{
	\refstepcounter{equation}
	\medskip
	\noindent\textbf{Examples~\theequation.} 
	\medskip
	\begin{enumerate}[(#1),ref=(#1),itemsep=5pt]
}
{
	\end{enumerate}\smallskip
}

\theoremstyle{remark}

\declaretheoremstyle[
spaceabove=\topsep, 
spacebelow=6pt,
headfont=\normalfont\itshape,
notefont=\normalfont, notebraces={(}{)},
bodyfont=\normalfont,
postheadspace=4pt,
qed=\mbox{\smaller[4]$\boxtimes$}
]{claimproofstyle}

\crefname{subsection}{Subsection}{Subsections}

\crefformat{footnote}{#2\footnotemark[#1]#3}

\crefrangeformat{enumi}{#3#1#4--#5#2#6}

\theoremstyle{plain}

\usepackage{color}

\definecolor{gris}{RGB}{90,90,90}

\definecolor{purple}{RGB}{116,0,159}

\renewcommand{\theequation}{\thesection.\arabic{equation}}

\makeatletter
\renewcommand\subsection{\@startsection{subsection}{2}%
	\z@{-1.5em}{.7em}%
	{\noindent\bfseries}}
\makeatother

\makeatletter

\def\l@section{\@tocline{1}{5pt}{0pc}{}{}}

\renewcommand{\tocsection}[3]{%
	\indentlabel{\@ifnotempty{#2}{\makebox[20pt][l]{%
				\ignorespaces#1 #2.\hfill}}}\sc #3\dotfill}

\newdimen{\tocsubsecmarg}
\def\l@subsection{\@tocline{2}{3pt}{0pc}{\tocsubsecmarg}{}}

\renewcommand{\tocsubsection}[3]{%
	\indentlabel{\@ifnotempty{#2}{\makebox[30pt][l]{%
				\ignorespaces#1 #2.\hfill}}}#3\dotfill}
\makeatother

\let\oldtocsubsection=\tocsubsection

\renewcommand{\tocsubsection}[2]{\hspace{3em} \oldtocsubsection{#1}{#2}}

\newcommand{\pB}[2][\BC]{#1^{\otimes #2}}

\newcommand{\Inv}{\text{\textnormal{Inv}}_a(X,\nu)}

\newcommand{\Lnorm}[1]{\norm{#1}_{L^2}}

\newcommand{\IPstar}{\text{\textnormal{IP}}^\ast}

\newcommand{\meanf}{\tilde{f}}

\title[]{Mixing and double recurrence in probability groups}
\author[]{Anush Tserunyan}
\thanks{The author was partially supported by NSF Grant DMS-1501036 and NSERC Discovery Grant RGPIN-2020-07120.}
\address{Department of Mathematics and Statistics, McGill University, Montreal, QC, Canada}
\email{anush.tserunyan@mcgill.ca}



\date{}

\begin{document}

\begin{abstract}
We define a class of groups equipped with an invariant probability measure, which includes all compact groups and is closed under taking ultraproducts with the induced Loeb measure; in fact, this class also contains the ultraproducts of all locally compact unimodular amenable groups. We call the members of this class \textit{probability groups} and develop the basics of the theory of their probability-measure-preserving actions, including a natural notion of mixing. A short proof reveals that for probability groups mixing implies double recurrence, which generalizes a theorem of Bergelson and Tao proved for ultraproducts of finite groups. Moreover, a quantitative version of our proof gives that $\e$-approximate mixing implies $3\sqrt{\e}$-approximate double recurrence. Examples of approximately mixing probability groups are quasirandom groups introduced by Gowers, so the last theorem generalizes and sharpens the corresponding results for quasirandom groups of Bergelson and Tao, as well as of Austin. Lastly, we point out that the fact that the ultraproduct of locally compact unimodular amenable groups is a probability group provides a general alternative to Furstenberg correspondence principle.
\end{abstract}


\
\vspace{-.7cm}

\maketitle


\tableofcontents

\vspace{-.5cm}

\begin{remarklike*}{Note to the reader}
	The reader who prefers to focus on finite or compact groups -- without going into the definition of general probability groups -- may safely skip the first two sections and read the rest having compact or finite groups in mind in lieu of probability groups.
		
\end{remarklike*}

\section{Overview of ultraproducts and Loeb measures}

We start with a quick overview of the construction of ultraproduct of measure spaces and discuss involved measurability issues, and thus motivate our definitions below, which otherwise might seem overly complicated.

\subsection{Ultraproducts}

For a short yet thorough survey of ultraproducts, we refer the reader to \cite{Keisler:survey:ultraproducts}.

Let $I$ be a countable index set and let $\alpha$ be an \emph{ultrafilter} on $I$, by which we mean a finitely additive $\set{0,1}$-valued measure defined on all of $\Pow(I)$. To make what follows nontrivial, we also assume that the ultrafilter $\alpha$ is \emph{nonprincipal}, i.e. is not a Dirac point measure (in particular, finite sets are $\alpha$-null). For a sequence $(X_i)_{i \in I}$ of sets, we think of elements $x,y$ of the product $\prod_{i \in I} X_i$ as functions $x,y : I \to \bigcup_{i \in I} X_i$, and thus, define the following equivalence relation
\[
x =_\alpha y \defequivlong x(i) = y(i) \text{ for $\alpha$-a.e. $i \in I$}
\]
just like we do with functions on a measure space. We call the quotient space $X \defeq \prod_{i \in I} X_i / =_\alpha$ the \emph{ultraproduct} of $(X_i)_{i \in I}$ over $\alpha$ and denote it by $\prod_{i \to \alpha} X_i$. Continuing the analogy with usual measurable functions, we identify $x \in \prod_{i \in I} X_i$ with its equivalence class $[x]_\alpha$; likewise, we often identify a subset $S$ of $\prod_{i \in I} X_i$ with the union $[S]_\alpha$ of the equivalence classes of the elements of $S$.

One can think of the ultraproduct as a limit of the sets $X_i$, and, as such, it inherits the properties and structure enjoyed by $\alpha$-a.e. $X_i$. For example, if each $X_i$ is actually a group $(G_i, e_i, \cdot_i)$, then so is their ultraproduct: simply define the multiplication coordinate-wise and $(e_i)_{i \in I}$ would be the identity. More generally, \Los's theorem \cite{Keisler:survey:ultraproducts}*{Theorem 3.1} states that this is true for any first-order property. Moreover, this is sometimes true for non-first-order properties such as being a probability space; that is, given that each $X_i$ admits a probability measure $\mu_i$, one can build a limit probability measure on the ultraproduct, called the \emph{Loeb measure}. To describe this construction, we first need to discuss the main property of ultraproducts, namely, \emph{countable compactness}.

\subsection{Countable compactness}

Call a set $B \subseteq X$ a \emph{quasibox} (more commonly called an \emph{internal set}) if it is of the form $[\prod_{i \in I} B_i]_\alpha$, where $B_i \subseteq X_i$. Note that the collection of quasiboxes is an algebra: indeed, the closure under finite intersections is obvious and, perhaps somewhat counterintuitively, the complement of $[\prod_{i \in I} B_i]_\alpha$ is $[\prod_{i \in I} B_i^c]_\alpha$. Thus, quasiboxes form a clopen basis for the topology they generate.

Assume further that $\alpha$ is \emph{nonprincipal}, i.e. not a point-measure. Then, we get the main property of ultraproducts, namely \emph{countable compactness} (also known as \emph{countable saturation}), which exhibits them as a certain kind of compactification.

\begin{prop}[Countable compactness]\label{ctbl-compactness}
	For any countable collection $\CC$ of quasiboxes in $X \defeq \prod_{i \to \alpha} X_i$, the topology on $X$ generated by $\CC$ is compact.
\end{prop}
\begin{proof}
	Let $\AC$ be the algebra generated by $\CC$ and note that $\AC$ is still countable and that it is enough to show that the topology generated by $\AC$ is compact. To show the latter, it is enough to prove that any sequence $(B^{(n)})_{n \in \N}$ of quasiboxes with the finite intersection property has nonempty intersection. Writing $B^{(n)} = [\prod_{i \in I} B^{(n)}_i]_\alpha$, we see that, for each $N \in \N$, for $\alpha$-a.e. $i \in I$,
	$
	\bigcap_{n < N} B^{(n)}_i \ne \0.
	$
	Identifying $I \defeq \N$, for each $i \in I$, let $N_i$ be the largest number $\le i$ such that $\bigcap_{n < N_i} B^{(n)}_i \ne \0$ and, using the axiom of choice, pick $x_i$ from $\bigcap_{n < N_i} B^{(n)}_i$. We claim that $x \defeq (x_i)_{i \in I}$ belongs to $B^{(N)}$, for every $N \in \N$. Indeed, because $\bigcap_{n \le N} B^{(n)} \ne \0$, we have that, for $\alpha$-a.e. $i \in I$, $\bigcap_{n \le N} B^{(n)}_i \ne \0$, and hence, $N_i \ge \min\set{N,i}$. Because $\alpha$ is nonprincipal, $i > N$ for $\alpha$-a.e. $i \in I$, so $N_i \ge N$, and hence, $x_i \in \bigcap_{n \le N} B^{(n)}_i$.
\end{proof}

\subsection{The Loeb measure construction}

A witty application of countable compactness is a construction of a countably additive measure on the ultraproduct of (even just finitely additive) measure spaces due to Loeb \cite{Loeb:measure}.

For each $i \in I$, let $(X_i, \BC_i, \mu_i)$ be a finitely additive measure space. Let $X \defeq \prod_{i \to \alpha} X_i$ and let $\AC \defeq \prod_{i \to \alpha} \BC_i$ denote the collection of all quasiboxes in $X$ with sides from the $\BC_i$, i.e. $[\prod_{i \in I} B_i]_\alpha$ with $B_i \in \BC_i$ for each $i \in I$. Clearly, $\AC$ is an algebra and the following defines a finitely additive measure on it:
\begin{equation}\label{eq:Lob-measure}
	\rho\left([\prod_{i \in I} B_i]_\alpha\right) \defeq \lim_{i \to \alpha} \mu_i(B_i).
\end{equation}
This limit is well-defined and it always exists because the space $[0,+\w]$ is compact.

Let $\BC \defeq \sigma(\AC)$ be the $\sigma$-algebra on $X$ generated by $\AC$; we refer to $\BC$ as the \emph{Loeb $\sigma$-algebra induced by the $\BC_i$}. We would like to extend $\mu$ to $\BC$ using the Caratheodory extension theorem. To do so, one only has to check that $\rho$ is countably additive on $\AC$, i.e. whenever a set $A \in \AC$ is a countable disjoint union of a sequence of nonempty sets $A_n \in \AC$, $n \in \N$, the measure $\rho(A)$ is equal to $\sum_{n \in \N} \mu(A_n)$. But this never occurs because the topology generated by $A^c$ and the sets $A_n$ is compact by \cref{ctbl-compactness}. Thus, we just proved the following.

\begin{prop}[Loeb]\label{Loeb-meas_ctbly_additive}
	The ultraproduct $X$ of finitely additive measure spaces $(X_i, \BC_i, \mu_i)$ admits a countably additive measure $\mu$ on the $\sigma$-algebra $\BC$ generated by the quasiboxes $[\prod_{i \in I} B_i]_\alpha$ with $B_i \in \BC_i$, on which $\mu$ is defined as in \labelcref{eq:Lob-measure}.
\end{prop}

We refer to this $\mu$ as the \emph{Loeb measure}.

\section{Probability groups and their actions}

The main goal of this section is to define a class of groups with an invariant probability measure, so that this class is closed under ultraproducts and contains all compact groups\refstepcounter{footnote}\footnote{Here and below by a compact group we mean a compact Hausdorff topological group.}. 

\subsection{Fubini systems}

Our global goal is to define a class of groups equipped with an invariant probability measure such that this class contains all compact groups and is closed under taking ultraproducts. Let us see what happens when we take the ultraproduct of finite groups; more precisely, for each $i \in I$, consider $(G_i, \BC_i, \mu_i)$, where $G_i$ is a finite group, $\BC_i \defeq \Pow(G_i)$, and $\mu_i$ the Haar measure (i.e.~normalized counting measure). We equip the ultraproduct $G$ of $(G_i)_{i \in \N}$ with a Loeb $\sigma$-algebra $\BC$ and the Loeb measure $\mu$ on $\BC$. As mentioned above, $G$ is also a group. However, we have an issue with measurability of the group operation on $G$.

\begin{notation}
	For a set $X$ and a $\sigma$-algebra $\BC$ on $X$, denote by $\pB{k}$ the $\sigma$-algebra on $X^k$ generated by $\BC^k$.
\end{notation}

Note that for each $i \in I$, the multiplication operation on $G_i$ is measurable as a function from $(G_i^2, \pB[\BC_i]{2})$ to $(G_i, \BC_i)$. However, the multiplication on $G$ need not be measurable as a function $(G^2, \pB{2}) \to (G, \BC)$. The reason is that $\pB{2}$ is in general a strictly smaller $\sigma$-algebra than the Loeb $\sigma$-algebra $\BC^{(2)}$ induced by the sequence $(\pB[\BC_i]{2})_{i \in I}$; the first example showing the strictness was given by Hoover \cite{Hoover:strictness_Loeb-alg_example} (see also \cite{Albeverio-at-al:book:nonstandard_methods}*{Example 3.2.13} for an exposition by D. Norman) and it was later shown in general for atomless probability spaces by Sun \cite{Sun:strictness_Loeb-alg_atomless}*{Proposition 6.6}. By \Los's theorem, the multiplication operation on $G$ is indeed $\BC^{(2)}$-measurable, and although $\BC^{(2)}$ is larger than $\pB{2}$, it is not that far from $\pB{2}$ in the sense that Fubini's theorem still holds, see \cite{Keisler:book:infinitesimal_stoch_analysis}*{1.14b} and \cite{Hurd-Loeb:book}*{Theorem 5.5}. The following definition captures this structure.

\begin{defn}\label{defn:Fubini_system}
	Let $X$ be a set. For each $k \ge 1$, let $\BC^{(k)}$ be a $\sigma$-algebra on $X^k$ and let $\mu^{(k)}$ be a (countably additive) probability measure on $\BC^{(k)}$. The tuple $(X, (\BC^{(k)}, \mu^{(k)})_{k \ge 1})$ is called a \emph{symmetric Fubini probability system} if, for each $k,l,n \ge 1$, the following conditions hold:
	\begin{enumref}{defn:Fubini_system}
		\item (symmetry) the coordinate permutation maps on $X^k$ are measurable and $\mu^{(k)}$-preserving;
		
		\item\label{defn:Fubini_system:projections} the \emph{projection} $(x, y) \mapsto x : X^{k+l} \to X^k$ is measurable and measure-preserving; equivalently, $\BC^{(k+l)} \supseteq \BC^{(k)} \times \BC^{(l)}$ and $\mu^{(k+l)} \rest{\BC^{(k)} \times \BC^{(k)}} = \mu^{(k)} \times \mu^{(l)}$;
		
		\item\label{defn:Fubini_system:duplications} the \emph{duplicating} map $(x_1, x_2,\hdots,x_k) \mapsto (x_1, x_1, x_2, \hdots, x_k) : X^k \to X^{k+1}$ is measurable;
		
		\item \label{defn:Fubini_system:Fubini} for every $A \in \BC^{(k+l)}$, the \emph{Fubini property} holds, namely:
		
		\begin{enumerate}[(a)]
			\item for every $x \in X^k$, the fiber $A_x$ is in $\BC^{(l)}$;

			\item the function $x \mapsto \mu^{(l)}(A_x) : X^k \to \R$ is measurable;
			
			\item \label{defn:Fubini_system:Fubini:integral} $\mu^{(k+l)}(A) = \displaystyle \int_{X^k} \mu^{(l)}(A_x) d\mu^{(k)}(x)$.
		\end{enumerate}
	\end{enumref}
\end{defn}

Similar definitions have been given in \cite{Keisler:prob_quantifiers}, \cite{Bagheri-Pourmahdian:logic_integration}, and \cite{Goldbring-Towsner:approx_logic_for_measures}.

\begin{obs}
	In the definition of Fubini systems, the symmetry of the $\sigma$-algebras implies that property \labelcref{defn:Fubini_system:duplications} holds for functions duplicating any $x_i$, not just $x_1$.
\end{obs}

For the sake of examples below, we also define more general \textit{\textbf{finitely additive} symmetric Fubini probability systems} the same way as in \cref{defn:Fubini_system} except that the measures $\mu^{(k)}$ are only finitely additive and the integral in \labelcref{defn:Fubini_system:Fubini}\labelcref{defn:Fubini_system:Fubini:integral} is understood as the unique mean\footnote{A \textit{mean} $\mu$ on a space of bounded functions (including the constant functions) is a linear functional such that $\mu(f) \ge 0$ if $f \ge 1$ and $\mu(1) = 1$.} on $L^\infty(X^k, B^{(k)}, \mu^{(k)})$ extending $\mu^{(k)}$, which we also denote by $\mu^{(k)}$ below.

\subsection{Probability groups}

\begin{defn}\label{defn:prob_group}
A (finitely additive) symmetric Fubini probability system $\big(G, (\BC^{(k)}, \mu^{(k)})_{k \ge 1}\big)$ is called a (\textit{finitely additive}) \emph{probability group} if $G$ is a group such that 
\begin{enumref}{defn:prob_group}
	\item\label{defn:prob_group:mult-inverse} for any $k \ge 1$, the left \emph{multiplication action} of $G$ on the first coordinate of $G^k$ and the \emph{inversion} of the first coordinate are measurable; more precisely, the maps
	\[
	(g_0,g_1,g_2,\hdots,g_k) \mapsto (g_0 g_1, g_2, \hdots, g_k) : G^{k+1} \to G^k
	\]
	and 
	\[
	(g_1,g_2,\hdots,g_k) \mapsto (g_1^{-1},g_2,\hdots,g_k) : G^k \to G^k
	\]
	are measurable;
	
	\item\label{defn:prob_group:measure-invariance} $\mu^{(1)}$ is \emph{invariant} under the two-sided multiplication and inverse; more precisely, for any $A \in \BC^{(1)}$, 
	\[
	\mu^{(1)}(g \cdot A) = \mu^{(1)}(A \cdot g) = \mu^{(1)}(A) = \mu^{(1)}(A^{-1}).
	\]
\end{enumref}
\end{defn}

\begin{remarklike}{Historical remark}
	The author was surprised to find a very similar definition in \cite{Weil} as it does not seem like Weil applies it to ultraproducts, which is where having a stronger $\sigma$-algebra on the product is needed.
\end{remarklike}

Below, we often simply write $G$ or $(G,\mu)$ for a probability group when the $\sigma$-algebras and the measures on higher dimensions are understood or not important for the discussion.

\begin{examples}\label{examples:prob_gps}
\item \textit{Every finite group is a probability group} (with the normalized counting measures).

\item \label{examples:prob_gps:compact} More generally, \textit{every compact Hausdorff group $G$ with is a probability group}. Here the $\sigma$-algebras $\BC^{(k)}$ are the Borel $\sigma$-algebra of the topology of $G^k$ and the measures $\mu^{(k)}$ are the unique normalized Haar measures on $\BC^{(k)}$.

\item \textit{Every countable amenable group $G$ is a finitely additive probability group}. Here the $\sigma$-algebras are just the powersets, but the measures take a bit to describe. 

By \cite{Greenleaf:book}*{Lemmas 1.1.1, 1.1.3}, there is a two-sided invariant mean $\mu^{(1)}$ on $\ell^{\infty}(G)$. This $\mu^{(1)}$ is a weak*-limit in $\ell^{\infty}(G)^\ast$ of a sequence $(\nu_n)$ of (finitely supported) probability measures $\nu_n$ on $G$. By the proofs of \cite{Greenleaf:book}*{Theorems 2.4.2, 2.4.3} and using the two-sided invariance of $\mu^{(1)}$, we may assume (passing to a subsequence) that for each $g \in G$, both $\norm{g \ast \nu_n - \nu_n}_1$ and $\norm{\nu_n \ast g - \nu_n}_1$ converge to $0$, where $g \ast \nu$ and $\nu \ast g$ are the pushforwards of $\nu$ by the maps on $G$ of left multiplication by $g^{-1}$ and right multiplication by $g$, respectively. (The proof of \cite{Greenleaf:book}*{Theorem 2.4.3} goes through for multiplication on both sides because the counting measure on $G$ is two-sided invariant.)

Now for each $k \ge 1$, there is a subsequence of $(\nu_n^k)$ that converges in the weak* topology of $\ell^\infty(G)^\ast$ and we can ensure (recursively) that these subsequences are nested. (Alternatively, by an Arzelà--Ascoli style diagonalization, we could get one subsequence that works for all $k$.) Let $\mu^{(k)}$ be the limit of the \onum{k} subsequence. The nestedness of the subsequences ensures the Fubini property \labelcref{defn:Fubini_system:Fubini} because this property holds for the $\nu_n^{(k)}$. Furthermore, $\norm{g \ast \nu_n^{(k)} - \nu_n^{(k)}}_1$ and $\norm{\nu_n^{(k)} \ast g - \nu_n^{(k)}}_1$ converge to $0$, which implies that $\mu^{(k)}$ is two-sided invariant.

\item More generally, \textit{every locally compact Hausdorff unimodular amenable group $G$ is a finitely additive probability group}. Here, for each $k \ge 1$, the $\sigma$-algebra $\BC^{(k)}$ is the Borel $\sigma$-algebra of the topology of $G^k$, and the measure $\mu^{(k)}$ is defined as in the previous example with the following modifications. 

By the last paragraph of \cite{Greenleaf:book}*{\S 2.2}, there is a two-sided invariant mean $\mu^{(1)}$ on $L^\infty(G, \chi)$, where $\chi$ is a two-sided invariant Haar measure on $G$. This $\mu^{(1)}$ is a weak*-limit in $\ell^{\infty}(G)^\ast$ of a net $(\nu_i)_{i \in I}$ of (compactly supported) probability measures $\nu_n$ on $G$. By the proofs of \cite{Greenleaf:book}*{Theorems 2.4.2, 2.4.3}, passing to a subnet, we may assume that for each $g \in G$, both $\norm{g \ast \nu_i - \nu_i}_1$ and $\norm{\nu_i \ast g - \nu_i}_1$ converge to $0$ (the proof of \cite{Greenleaf:book}*{Theorem 2.4.3} goes through for multiplication on both sides because $\chi$ is two-sided invariant). The rest of the construction is the same as in the previous example, with subsequences replaced with subnets.
\end{examples}

The following statement is the main reason for defining probability groups as it provides a plethora of important examples and a tool for proving statements about (finitely additive) probability groups.

\begin{prop}\label{closure_under_ultraproducts}
The class of probability groups is closed under ultraproducts. In fact, the ultraproduct of finitely additive probability groups equipped with the induced Loeb measures is a (countably additive) probability group.
\end{prop}
\begin{proof}
Let $\alpha$ be a nonprincipal ultrafilter on $I \defeq \N$ and for each $i \in I$, let $\big(G_i, (\BC_i^{(k)}, \mu_i^{(k)})_{k \ge 1}\big)$ be a finitely additive probability group. Take $G \defeq \prod_{i \to \alpha} G_i$ and for each $k \ge 1$, $\BC^{(k)} \defeq \sigma\left(\prod_{i \to \alpha} \BC_i^{(k)}\right)$ and $\mu^{(k)} \defeq \prod_{i \to \alpha} \mu_i^{(k)}$. It is now not hard to verify that we have obtained a probability group. Indeed, it is a theorem of Keisler \cite{Keisler:book:infinitesimal_stoch_analysis}*{1.14b} and of Hurd and Loeb \cite{Hurd-Loeb:book}*{Theorem 5.5} that the Fubini property holds. Checking the rest of the conditions of \cref{defn:prob_group} amounts to straightforward applications of \Los's theorem. Finally, the countable additivity of the Loeb measure is given by \cref{Loeb-meas_ctbly_additive}.
\end{proof}

Combined with \multiexamplesref{examples:prob_gps}, we now get the following examples, in which the $\sigma$-algebras and measures on the ultraproduct are given by the Loeb measure construction.

\begin{examples}\label{examples:ultra_prob_gps}

\item Every ultraproduct of compact Hausdorff groups is a probability group.

\item \label{examples:ultra_prob_gps:ctbl-amenable} Every ultraproduct of countable amenable groups is a probability group.

\item \label{examples:ultra_prob_gps:lc-amenable} More generally, every ultraproduct of locally compact Hausdorff unimodular amenable groups is a probability group.\footnote{Thanks to the anonymous referee for suggesting this example.} 

\end{examples}

\examplesref{examples:ultra_prob_gps}{examples:ultra_prob_gps:ctbl-amenable} can be viewed as an alternative to the \textit{Furstenberg correspondence principle}; indeed, \Los's theorem transfers back-and-forth first-order statements about a countable amenable group and its ultrapower. However, the ultrapower has the advantage of having its measure be \textbf{countably additive}. Similarly, \examplesref{examples:ultra_prob_gps}{examples:ultra_prob_gps:lc-amenable} can be viewed as a more general alternative to Furstenberg correspondence that applies to all locally compact unimodular amenable groups.

\subsection{Properties of probability groups}

\begin{notation}
	For $G$ a group, define its $i^{\text{th}}$ coordinate left and right actions on $G^k$ by $g \cdot_i^k (g_1,\hdots,g_i,\hdots,g_k) \defeq (g_1,\hdots,g g_i,\hdots,g_k)$ and $(g_1,\hdots,g_i,\hdots,g_k) \cdot_i^k g \defeq (g_1,\hdots, g_i g,\hdots,g_k)$; denote the action functions by $L_i^k : G^{k+1} \to G^k$ and $R_i^k : G^{k+1} \to G^k$. Similarly, define the $i^{\text{th}}$ coordinate inverse action on $G^k$ by $I_i^k : (g_1,\hdots,g_i,\hdots,g_k) \mapsto (g_1,\hdots,g_i^{-1},\hdots,g_k)$.
\end{notation}

\begin{obs}
	In a probability group $\big(G, (\BC^{(k)}, \mu^{(k)})_{k \ge 1}\big)$, because the $\BC_k$ are symmetric, it follows that for every $k \ge 1$ and $i \le k$, the maps $L_i^k, R_i^k : (G^{k+1}, \BC^{(k+1)}) \to (G^k, \BC^{(k)})$ and $I_i^k : (G^k, \BC^{(k)}) \to (G^k, \BC^{(k)})$ are measurable.
\end{obs}

\begin{prop}[Invariance in all dimensions]\label{invariance_in_all_dim}
	In a probability group $\big(G, (\BC^{(k)}, \mu^{(k)})_{k \ge 1}\big)$, for every $k \ge 1$, the measure $\mu^{(k)}$ is invariant under the left/right multiplication and inverse actions on any coordinate, i.e. for any $A \in \BC^{(k)}$, $i \le k$, and $g \in G$,
	\[
	\mu^{(k)}(g \cdot_i^k A) = \mu^{(k)}(A \cdot_i^k g) = \mu^{(k)}(I_i^k(A)) = \mu^{(k)}(A).
	\]
\end{prop}
\begin{proof}
	This is due to the Fubini property. For example, because the function $L_1^k$ is measurable, its fiber $(L_1^k)_g$ is also measurable for any fixed $g \in G$, which implies that for any $A \in \BC^{(k)}$, $g \cdot_1^k A \in \BC^{(k)}$. Moreover, by the Fubini property and the invariance of $\mu^{(1)}$ under the action of $G$, putting $h \defeq (g_2, \hdots, g_k)$, we have
	\begin{align*}
		\mu^{(k)}(g \cdot_1^k A) 
		&= 
		\int_{G^{k-1}} \mu^{(1)}\big((g \cdot_1^k A)_h\big) d\mu^{(k-1)}(h) 
		\\ 
		&= 
		\int_{G^{k-1}} \mu^{(1)}(g \cdot A_h) d\mu^{(k-1)}(h) 
		\\
		&= 
		\int_{G^{k-1}} \mu^{(1)}(A_h) d\mu^{(k-1)}(h) 
		= 
		\mu^{(k)}(A). \qedhere
	\end{align*}
\end{proof}

\begin{prop}[Word maps]\label{word-maps}
	In any probability group $(G, (\BC^{(k)}, \mu^{(k)})_{k \ge 1})$, all word multiplication maps are measurable; more precisely, for any $n,k \ge 1$ and any words $w_1,w_2,\hdots,w_k$ in the alphabet 
	\[
	\Sigma \defeq \set{x_1,x_2,\hdots,x_n,y_1,y_2,\hdots,y_n},
	\]
	the map 
	\[
	\vec{g} \mapsto \big(w_1(\vec{g}), w_2(\vec{g}),\hdots,w_k(\vec{g})\big)
	:
	G^n \to G^k
	\]
	is measurable, where, for a word $w \in \Fin[\Sigma]$, $w(\vec{g})$ is the result of plugging in $x_i \defeq g_i$, $y_i \defeq g_i^{-1}$ in $w$ and multiplying out.
\end{prop}
\begin{proof}
	Instead of giving a notation-heavy proof for the general case, we do it for the map $(u,x,y,z) \mapsto (y^2 z x^{-1}, u^{-1} x^2) : G^4 \to G^2$. The map
	\[
		(u,x,y,z) \mapsto (y, y, z, x, u, x, x)
	\]
	is measurable due to iterative applications of \labelcref{defn:Fubini_system:duplications} and symmetry. Similarly, \labelcref{defn:prob_group:mult-inverse} implies that the maps
	\[
	(y, y, z, x, u, x, x) \mapsto (y, y, z, x^{-1}, u^{-1}, x, x) \mapsto (y^2 z x^{-1}, u^{-1} x^2)
	\]
	are measurable, so taking their composition finishes the proof.
\end{proof}

\subsection{Measure-preserving actions of probability groups}

We will now define a natural class of actions for probability groups. We again have a measurability issue to deal with, which makes the definition very similar to the definitions of Fubini systems and probability groups put together. Thus, we will give a rather informal definition instead, hoping that the suppressed details are understood.

\begin{defn}\label{defn:actions_of_prob_gps}
Let $(G, (\BC^{(k)}, \mu^{(k)})_{k \ge 1})$ be a probability group, $(X, \CC, \nu)$ a probability space, and let $(x,g) \mapsto x \cdot_a g : X \times G \to X$ be a right action of $G$ on $X$. We call this action \emph{measure-preserving} with respect to $\sigma$-algebras $\CC^{(k)}$ on $X \times G^{k-1}$ with $\CC^{(1)} \defeq \CC$ and probability measures $\nu^{(k)}$ on $\CC^{(k)}$ with $\nu^{(1)} \defeq \nu$, $k \ge 1$, such that 
\begin{enumref}{defn:actions_of_prob_gps}
\item the natural extensions\footnote{We mean that these maps leave the $X$-coordinate unchanged.} from $G^k$ to $X \times G^k$ of all of the \emph{permutation}, \emph{projection}, \emph{duplicating}, \emph{group multiplication} and \emph{inversion} maps are measurable with respect to the corresponding $\CC^{(k)}$-s and the permutation and projection maps are \emph{measure-preserving}; in particular, $\CC^{(k+l)} \supseteq \CC^{(l)} \otimes \BC^{(k)}$ and $\nu^{(k+l)} \rest{\CC^{(l)} \otimes \BC^{(k)}} = \nu^{(l)} \times \mu^{(k)}$;

\item the maps $(x, g_1, g_2, \hdots, g_k) \mapsto (x \cdot_a g_1, g_2, \hdots, g_k) : X \times G^k \to X \times G^{k-1}$ is measurable;

\item the action \emph{preserves the measure} $\nu$, i.e. $\nu(A \cdot_a g^{-1}) = \nu(A)$ for all $g \in G$ and $A \in \CC$;

\item the \emph{Fubini property} holds in all dimensions.
\end{enumref}

We do not mention the $\sigma$-algebras $\CC^{(k)}$ and the measures $(\mu^{(k)})$ if it is not important for the discussion, but a measure-preserving action, by definition, comes equipped with this data. We also often simply write $G \actson (X,\nu)$ or $(G,\mu) \actson (X, \nu)$ for a measure-preserving action of a probability group $(G,\mu)$ on a probability space $(X,\nu)$.
\end{defn}

\begin{example}
For a probability group $G$, the left and right translation actions $x \cdot_\ell g \mapsto g^{-1} x$ and $x \cdot_r g \mapsto x g$, as well as the conjugation action $x \cdot_c g \mapsto g^{-1} x g = g \cdot_\ell (g \cdot_r x)$ of $G$ on itself, are measure-preserving (right) actions with respect to $\CC^{(k)} \defeq \BC^{(k)}$ and $\nu^{(k)} \defeq \mu^{(k)}$.
\end{example}

It is routine to verify that the natural analogues of \cref{invariance_in_all_dim,word-maps}, as well as the first part of \cref{closure_under_ultraproducts} (closedness under ultraproducts), hold for measure-preserving actions of probability groups on probability spaces.

\begin{defn}[Unitary representations]
A right action $a : G \actson (X,\nu)$ of an (abstract) group $G$ on a probability space $(X,\nu)$ by measure-preserving automorphisms induces a left action $G \actson L^2(X, \nu)$ by unitary operators. This action, still denoted by $\cdot_a$, is known as the \textit{Koopman representation} of the original action and is defined by
$
(g \cdot_a f)(x) = f(x \cdot_a g).
$
Let $\Inv \subseteq L^2(X,\nu)$ denote the subspace of functions $f$ invariant under this action, i.e. $g \cdot_a f = f$ for all $g \in G$. Finally, let
$
P_a : L^2(X, \nu) \to \Inv
$
be the orthogonal projection onto $\Inv$.
\end{defn}

Below we use $\gen{\cdot,\cdot}_X$ to denote the inner product in $L^2(X, \nu)$. All $L^2$-spaces and more generally, all Hilbert spaces are assumed to be complex.

\section{Ergodicity and mixing}

\subsection{The mean ergodic theorem}

\begin{defn}
	A measure-preserving action $a : (G,\mu) \actson (X,\nu)$ of a probability group on a probability space is called \emph{ergodic} if any measurable $a$-invariant subset of $X$ is either $\nu$-null or $\nu$-conull.
\end{defn}

If $G$ is a probability group and the action $a : G \actson G$ is either the left or right translation, then for $f \in L^2(G)$, $P_a(f)$ is just the mean of $f$ because these actions are transitive, so the only invariant functions are constants. In general, the following gives an explicit computation of $P_a$ for arbitrary measure-preserving actions of probability groups.

\begin{prop}[Mean ergodic theorem for probability groups]\label{mean_ergodic_theorem}
Let $a : (G,\mu) \actson (X,\nu)$ be a measure-preserving action of a probability group on a probability space. For all $f \in L^2(X,\nu)$,
\[
P_a(f)(x) = \int_G (g \cdot_a f)(x) d\mu(g).
\]
In particular, if the action is ergodic, then for $\nu$-a.e. $x \in X$,
\[
\int_G (g \cdot_a f)(x) d\mu(g) = \int_X f(y) d\nu(y).
\]
\end{prop}
\begin{proof}
Putting $\meanf(x) \defeq \int_G (g \cdot_a f)(x) d\mu(g)$ and fixing $\phi \in \Inv$, we need to show that $f - \meanf$ and $\phi$ are orthogonal, for which it is enough to show that $\gen{f, \phi}_X = \gen{\meanf, \phi}_X$. Compute:
\begin{align*}
\gen{\meanf, \phi}_X 
&= 
\int_X \int_G (g \cdot_a f)(x) \phi(x) d\mu(g) d\nu(x) 
\\
\eqcomment{Fubini} 
&= 
\int_G \gen{g \cdot_a f, \phi}_X d\mu(g) 
\\
\eqcomment{unitarity} 
&= 
\int_G \gen{f, g^{-1} \cdot_a \phi}_X d\mu(g) 
\\
\eqcomment{invariance of $\phi$} 
&= 
\int_G \gen{f, \phi}_X d\mu(g) 
= 
\gen{f, \phi}_X.
\end{align*}
Furthermore, if the action is ergodic, then the only functions in $\Inv$ are constants, so $\meanf \equiv \int_X f(x) d\nu(x)$ $\nu$-a.e.
\end{proof}

\subsection{Mixing}

For a measure $\mu$, we write $\forall^\mu$ to mean ``for $\mu$-a.e.''.

\begin{defn}
For a probability group $(G, \mu)$ and a probability space $(X, \nu)$, call a measure-preserving action $a : G \actson X$ \emph{mixing along $\mu$} (or just \emph{mixing}) if for any $f_1, f_2 \in L^2(X,\nu)$,
\[
(\forall^\mu g \in G) \ \gen{f_1, g \cdot_a f_2}_X = \gen{P_a(f_1), P_a(f_2)}_X.
\]
\end{defn}

One could also give an abstract definition of \emph{mixing along a filter $\FC \subseteq \Pow(G)$} for any group $G$ as follows: for any $f_1, f_2 \in L^2(X,\nu)$,
\[
\lim_{g \to \FC} \gen{f_1, g \cdot_a f_2}_X = \gen{P_a(f_1), P_a(f_2)}_X.
\]
For ergodic actions, this generalizes the usual notions of mixing such as
\begin{itemize}
\item \emph{weak mixing} for amenable $G$ with the filter $\FC$ of density-one sets;

\item \emph{mild mixing} for arbitrary discrete $G$ with filter $\FC \defeq \IPstar$;

\item \emph{strong mixing} for arbitrary discrete $G$ with the Fr\'{e}chet filter $\FC$.
\end{itemize}

In our case, due to the countable additivity of $\mu$, the definition of $\mu$-mixing is equivalent to mixing along the filter of $\mu$-conull sets. 

\begin{remark}
A similar definition of mixing along a filter for ergodic actions was considered by Tucker-Drob \cite[Chapter 7]{TuckerDrob:thesis}.
\end{remark}

\begin{example}[Ultra quasirandom groups]\label{example_ultraquasirandom_right_mixing}
In \cite{Bergelson-Tao:mult_rec_quasirandom}, the authors consider finite groups that are approximately mixing (i.e. mixing with a small error); more precisely, they consider so-called \emph{$D$-quasirandom groups}, introduced by Gowers in \cite{Gowers:quasirandom_gps}, that is: finite (or, more generally, compact Hausdorff) groups that do not admit any nontrivial unitary representations of dimension less than $D$ (\cref{def:quasirandom-groups} below). It is then shown that the right translation action of these groups on themselves is mixing with an error $D^{-1/2}$, with respect to the normalized Haar measure (see \cite{Bergelson-Tao:mult_rec_quasirandom}*{Proposition 3} or \cref{sec:quantitative_version} below). Therefore, taking an appropriate ultraproduct washes the error away, yielding a probability group whose right translation action on itself is genuinely mixing. More precisely, in \cite{Bergelson-Tao:mult_rec_quasirandom}, the authors define \emph{ultra quasirandom groups} as an ultraproduct of a sequence $(G_i, \mu_i)_{i \in \N}$ of finite groups, where $\mu_i$ is the normalized counting measure, each $G_i$ is $D_i$-quasirandom, and $D_i \to \w$. This is a probability group with respect to the induced Loeb measure, and, by \cite{Bergelson-Tao:mult_rec_quasirandom}*{Lemma 33}, its right translation action on itself is mixing.
\end{example}

We are finally ready to give the main definition, which at a glance may seem hard to check and unlikely to occur, but \cref{right_mixing_implies_mixing} below will settle the matter.

\begin{defn}
We call a probability group \emph{mixing} if all of its measure-preserving actions on probability spaces are mixing.
\end{defn}

\begin{prop}\label{right_mixing_implies_mixing}
A probability group $(G,\mu)$ is mixing if and only if its right translation action on itself is mixing.
\end{prop}
\begin{proof}
We show the nontrivial direction: suppose the right translation action $r : G \actson G$ is mixing and consider a measure-preserving action $a : G \actson X$ on a probability space $(X, \nu)$.

The idea is to switch from averaging over the action $a : G \actson X$ to averaging over the right translation action $r : G \actson G$; this is done using the Fubini property and the associativity of the action: for $g,h \in G$ and $x \in X$,
\[
	(x \cdot_a h) \cdot_a g =  x \cdot_a (h \cdot_r g).
\]

Turning to the actual proof, for a function $f : X \to \C$ and $x \in X$, let $f^{(x)} : G \to \C$ be defined by $g \mapsto (g \cdot_a f)(x)$. Observe that, for $g,h \in G$,
\begin{equation}\label{eq:fiber_func_assoc}
	\big(h \cdot_a (g \cdot_a f)\big)(x) = \big((hg) \cdot_a f\big)(x) = f^{(x)}(hg) = (g \cdot_r f^{(x)})(h).
\end{equation}
Fixing $f_1, f_2 \in L^2(X,\nu)$ and $g \in G$, we compute:
\begin{align*}
\gen{f_1, g \cdot_a f_2}_X 
&= 
\int_G \gen{f_1, g \cdot_a f_2}_X \, d\mu(h) 
\\
\eqcomment{unitarity} 
&= 
\int_G \gen{h \cdot_a f_1, h \cdot_a (g \cdot_a f_2)}_X \, d\mu(h)
\\
\eqcomment{by \labelcref{eq:fiber_func_assoc}} 
&=
\int_G \int_X f_1^{(x)}(h) \, (g \cdot_r f^{(x)})(h) \, d\nu(x) d\mu(h)
\\
\eqcomment{Fubini}
&= 
\int_X \gen{f_1^{(x)}, g \cdot_r f_2^{(x)}}_G \, d\nu(x).
\end{align*}
Because the right translation action is mixing and ergodic, we have
\[
(\forall x \in X) (\forall^\mu g \in G) \ \gen{f_1^{(x)}, g \cdot_r f_2^{(x)}}_G = (\int_G f_1^{(x)} d\mu) (\int_G f_2^{(x)} d\mu),
\]
so the Fubini property implies
\[
(\forall^\mu g \in G) (\forall^\nu x \in X) \ \gen{f_1^{(x)}, g \cdot_r f_2^{(x)}}_G = (\int_G f_1^{(x)} d\mu) (\int_G f_2^{(x)} d\mu).
\]
Moreover, the mean ergodic theorem (\cref{mean_ergodic_theorem}) applied to any $f \in L^2(X,\nu)$ gives $\int_G f^{(x)} d\mu = P_a(f)(x)$ for $\nu$-a.e. $x \in X$, so, for $\mu$-a.e. $g \in G$,
\[
\gen{f_1, g \cdot_a f_2}_X 
= 
\int_X P_a(f_1)(x) \, P_a(f_2)(x) \, d\mu(x) 
= 
\gen{P_a(f_1), P_a(f_2)}_X.
\qedhere
\]
\end{proof}

\begin{example}
As mentioned in \cref{example_ultraquasirandom_right_mixing}, the right translation action of an ultra quasirandom group on itself is mixing. Thus, ultra quasirandom groups are mixing probability groups. This, in particular, implies \cite{Bergelson-Tao:mult_rec_quasirandom}*{Lemma 34}.
\end{example}

\section{Double recurrence}

\begin{defn}
	Call a probability group $(G,\mu)$ \emph{doubly recurrent} if for any $f_1, f_2, f_3 \in L^\w(G,\mu)$,
	\begin{equation}\label{eq:dbl_rec}
		(\forall^\mu g \in G) \ \int_G f_1 (g \cdot_\ell f_2) (g \cdot_c f_3) d\mu = \int_G f_1 P_\ell(f_2) P_c(f_3) d\mu,
	\end{equation}
	where $\cdot_\ell$ and $\cdot_c$ are, respectively, the left translation and the conjugation actions of $G$ on itself.
\end{defn}

\subsection{Mixing implies double recurrence}

The following theorem is the main result of the paper. It generalizes \cite{Bergelson-Tao:mult_rec_quasirandom}*{Theorem 41} proven for ultra quasirandom groups.

\begin{theorem}\label{dbl_rec}
Every mixing probability group is doubly recurrent.
\end{theorem}

Using transfer principle (or equivalently, considering an ultraproduct of counterexample quasirandom groups with $D \to \w$), Bergelson and Tao show in \cite{Bergelson-Tao:mult_rec_quasirandom}*{Theorem 5} that this theorem for ultra quasirandom groups implies approximate double recurrence for finite quasirandom groups with an implicit bound on the error. \cite{Bergelson-Tao:mult_rec_quasirandom}*{Corollary 7} interprets this in terms of the distribution of the quadruples $(g,x,gx,xg)$ with $x,g$ drawn uniformly and independently at random. See also \cite{Bergelson-Tao:mult_rec_quasirandom}*{Corollary 8} for a density noncommutative Schur theorem for quasirandom groups.

\medskip

Before going into the proof, we briefly explain its idea.

\begin{remarklike}{Idea of proof}
If we remove one of the factors $f_1$, $g \cdot_\ell f_2$ or $g \cdot_c f_3$ from \labelcref{eq:dbl_rec}, i.e. ``drop the degree'' of the product, then the equality would easily follow from single recurrence, i.e. mixing. We get rid of the factor $f_1$ and here is how. Linearity reduces to the orthogonal cases $P_c(f_3) = f_3$ and $P_c(f_3) = 0$, and the proof of the former case follows from the left translation action being mixing, so we are left with the case $P_c(f_3) = 0$. Assuming this, what we need to show is
\[
\forall^\mu g \ \gen{f_1, e_g}_G = 0,
\]
where $e_g = (g \cdot_\ell f_2) (g \cdot_c f_3)$. But the latter would follow basically from Bessel's inequality if we could show that $\set{e_g}_{g \in G}$ is an a.e.-orthogonal family in $L^2(G,\mu)$, i.e. 
\[
\forall^{\mu_2}(g,h) \gen{e_g, e_h}_G = 0.
\]
By the Fubini property and a change of variable, this is equivalent to
\[
\forall^\mu h \forall^\mu g \gen{e_g, e_{gh}}_G = 0,
\]
which, due to some regrouping and cancellation, easily follows from the right translation and the conjugation actions being mixing. This latter trick of replacing pairs $(g,h)$ by $(g,gh)$  is known as the \emph{van der Corput difference trick}, which can be thought of as an analogue of differentiation in this context because an application of this trick ``drops the degree''.
\end{remarklike}

\begin{remark}
In the proof of this theorem for an ultra quasirandom group given in \cite{Bergelson-Tao:mult_rec_quasirandom}, the authors restrict to a countable subgroup $\Gamma$ of $G$ and use an idempotent ultrafilter on $\Gamma$ as their notion of largeness, which is almost invariant under the translation action of $\Gamma$ on itself. We instead use the measure $\mu$ on $G$, or equivalently, the filter of $\mu$-conull sets, which is genuinely invariant and also has the advantage of being countably additive; the latter enables cleaner pigeon-hole arguments and replaces various limits with ``a.e.'' statements. The only price we pay is that our filter of $\mu$-conull sets is not ``ultra'', but this is not an issue as we can be careful enough to stay in the $\sigma$-algebra of measurable sets when needed.
\end{remark}

\subsection{Proof of \cref{dbl_rec}}

We start by recording a (cheap) Ramsey theorem for filters. For a filter $\FC$ on a set $X$, we write $\forall^\FC$ below to mean ``for an $\FC$-large set of points in $X$''.

\begin{lemma}[Ramsey for filters]\label{Ramsey}
Let $X$ be a set and $\FC$ a nonprincipal filter on it. If a set $R \subseteq X^2$ is such that
$
(\forall^\FC x \in X) \ (\forall^\FC y \in X) \ x R y,
$
then there is an infinite set $\set{x_n}_{n \in \N} \subseteq X$ such that $x_n R x_m$ for all $n < m$.
\end{lemma}
\begin{proof}
For each $x \in X$, let $R_x \defeq \set{y \in X : x R y}$. By the hypothesis, the set $A \defeq \set{x \in X : R_x \text{ is $\FC$-large}}$ is $\FC$-large. Put $A_0 = A$ and take $x_0 \in A_0$. Put $A_1 = R_{x_0} \cap A_0$ and note that $A_1$ is still $\FC$-large. Take $x_1 \in A_1$ distinct from $x_0$ (can do this because $\FC$ is nonprincipal). Repeat: put $A_2 = R_{x_1} \cap A_1$ and note that $A_2$ is still $\FC$-large. Take $x_2 \in A_2$ distinct from $x_0, x_1$; etc.
\end{proof}

We also recall the following basic Hilbert space fact, which follows from Bessel's inequality:

\begin{lemma}[Bessel]\label{Bessel}
Let $(e_n)_{n \in \N}$ be a bounded sequence of vectors in a Hilbert space $\Hil$. If the vectors in $(e_n)_{n \in \N}$ are pairwise orthogonal, then $\displaystyle\lim_{n \to \w} e_n = 0$ in the weak topology of $\Hil$, i.e. for every $f \in \Hil$, $\displaystyle\lim_{n \to \w} \gen{f, e_n} = 0$.
\end{lemma}

Putting this together with the Ramsey lemma applied to the filter of conull sets, we get a natural analogue of Bessel's lemma for measure:
\begin{lemma}[Random Bessel]\label{random_Bessel}
Let $(X, \mu)$ be a measure space with nonatomic $\mu \ne 0$ and let $(e_x)_{x \in X}$ be a bounded sequence\footnote{We use the term \textit{sequence} even when the index set is not $\N$.} in a Hilbert space $\Hil$. If
\[
(\forall^\mu x \in X) \ (\forall^\mu y \in X) \ \gen{e_x, e_y} = 0,
\]
then for every $f \in \Hil$,
$
(\forall^\mu x \in X) \ \gen{f, e_x} = 0.
$
\end{lemma}
\begin{proof}
Fix $f \in \Hil$ and suppose that the conclusion fails for this $f$. Then, there is $\e > 0 $ such that the set $Y = \set{x \in X : |\gen{f, e_x}| \ge \e}$ is not $\mu$-null (caution: $Y$ may not be measurable). Thus, the restriction of the filter of $\mu$-conull sets to $Y$ gives a nonprincipal filter $\FC$ on $Y$. Applying \cref{Ramsey} to $Y$ with filter $\FC$ and $R = \set{(x,y) \in Y^2 : \gen{e_x, e_y} = 0}$, we get an infinite bounded sequence $(e_{x_n})_{n \in \N}$ of pairwise orthogonal vectors such that for every $n \in \N$, $|\gen{f, e_{x_n}}| \ge \e$, contradicting \cref{Bessel}.
\end{proof}

Inviting group structure and Fubini to this party of Ramsey and Bessel, we get:

\begin{lemma}[Random van der Corput]\label{random_van_der_Corput}
Let $\big(G, (\BC^{(k)}, \mu^{(k)})_{k \ge 1}\big)$ be an infinite probability group and let $(e_g)_{g \in G}$ be a bounded sequence in a Hilbert space $\Hil$ such that the function $(g,h) \mapsto \gen{e_g, e_h} : G^2 \to \C$ is $\BC^{(2)}$-measurable. If 
\[
(\forall^\mu h \in G) \ (\forall^\mu g \in G) \ \gen{e_g, e_{gh}} = 0,
\]
then for all $f \in \Hil$,
$
(\forall^\mu g \in G) \ \gen{f, e_g} = 0.
$
\end{lemma}
\begin{proof}
By the Fubini property, $(\forall^\mu g \in G) \, (\forall^\mu h \in G) \ \gen{e_g, e_{gh}} = 0$. The invariance of $\mu$ allows for a change of variable $h \mapsto g^{-1} h$, yielding $(\forall^\mu g \in G) \, (\forall^\mu h \in G) \ \gen{e_g, e_h} = 0$, so the desired conclusion follows from \cref{random_Bessel}.
\end{proof}

\begin{remark}
This lemma has several cousins in the countable setting; e.g. for the filter on $\N$ of sets of density $1$ \cite{Furstenberg:book}*{Lemma 4.9}, for the filter on $\N$ of sets that meet every IP-set \cite{Furstenberg:book}*{Lemma 9.24} and for idempotent ultrafilters on countable groups \cite{Bergelson-McCutcheon:central_sets_Roth}*{Theorem 2.3}. A generalization of all of these statements is proven in \cite{me:vdC}*{Theorem 6.1}. See also \cref{approx-van-der-Corput} below for a quantitative version.
\end{remark}

We are now ready to prove the double recurrence theorem.

\begin{proof}[Proof of \upshape\cref{dbl_rec}]
Let $\big(G, (\BC^{(k)}, \mu^{(k)})_{k \ge 1}\big)$ be a mixing probability group. As we solely work in $G$, we omit the subscript $G$ from $\gen{\cdot,\cdot}_G$.

Because $g \cdot_c P_c(f_3) = P_c(f_3)$,
\[
\gen{f_1 (g \cdot_\ell f_2), g \cdot_c f_3} = \gen{f_1 (g \cdot_\ell f_2), g \cdot_c (f_3 - P_c(f_3))} + \gen{f_1 (g \cdot_\ell f_2), P_c(f_3)},
\]
so it is enough to prove the theorem in the following two orthogonal cases:

\begin{case}{1}{$P_c(f_3) = f_3$}
The desired identity \labelcref{eq:dbl_rec} turns into
\[
(\forall^\mu g \in G) \ \gen{f_1 f_3, g \cdot_\ell f_2} = \gen{f_1 f_3, P_\ell(f_2)},
\]
which immediately follows from the fact that the left translation action is mixing.
\end{case}

\begin{case}{2}{$P_c(f_3) = 0$}
Now identity \labelcref{eq:dbl_rec} turns into
\[
(\forall^\mu g \in G) \ \gen{f_1, (g \cdot_\ell f_2) (g \cdot_c f_3)} = 0,
\]
so it will follow from the random van der Corput lemma (\cref{random_van_der_Corput}) for $e_g = (g \cdot_\ell f_2) (g \cdot_c f_3)$ once we verify its hypothesis. It follows from the definition of probability groups (\cref{defn:prob_group}) that the function $G^2 \to \C$ defined by
\[
(g,h) \mapsto \gen{e_g, e_h} = \int_G (g \cdot_\ell f_2) \, (g \cdot_c f_3) \, (h \cdot_\ell f_2) \, (h \cdot_c f_3) \, d\mu
\]
is $\BC^{(2)}$-measurable. Furthermore, the sequence $(e_g)_{g \in G}$ in $L^2(G,\mu)$ is bounded because $f_2,f_3 \in L^\w(G,\mu)$ and $\mu$ is finite. It remains to verify that $\forall^\mu h \, \forall^\mu g \ \gen{e_g, e_{gh}} = 0$. To this end, we fix $h,g \in G$ and compute:
\begin{align*}
\gen{e_g, e_{gh}} 
&= 
\int_G (g \cdot_\ell f_2) \, (g \cdot_c f_3) \, ((gh) \cdot_\ell f_2) \, ((gh) \cdot_c f_3) \, d\mu 
\\
\eqcomment{associativity of actions and regrouping} 
&= 
\genbig{(g \cdot_\ell f_2) \, (g \cdot_\ell h \cdot_\ell f_2), (g \cdot_c f_3) \, (g \cdot_c h \cdot_c f_3)}
\\
\eqcomment{distributivity of actions over product} 
&= 
\genbig{g \cdot_\ell (f_2 (h \cdot_\ell f_2)), g \cdot_c (f_3 (h \cdot_c f_3))}
\\
\eqcommentlr{$\begin{array}{l}
F_2^{(h)} \defeq f_2 (h \cdot_\ell f_2) 
\\ 
F_3^{(h)} \defeq f_3 (h \cdot_c f_3)
\end{array}$}
&= 
\gen{g \cdot_\ell F_2^{(h)}, g \cdot_c F_3^{(h)}}
\\
\eqcomment{$g \cdot_c f = g \cdot_\ell g \cdot_r f$}
&= 
\gen{g \cdot_\ell F_2^{(h)}, g \cdot_\ell g \cdot_r F_3^{(h)}}
\\
\eqcomment{unitarity} 
&= 
\gen{F_2^{(h)}, g \cdot_r F_3^{(h)}}.
\end{align*}
Because the right translation action is mixing, we have that for every $h \in G$:
\[
(\forall^\mu g) \ \gen{F_2^{(h)}, g \cdot_r F_3^{(h)}} = (\int_G F_2^{(h)} d\mu) (\int_G F_3^{(h)} d\mu).
\]
But the conjugation action is mixing as well, so
\[
(\forall^\mu h) \ \int_G F_3^{(h)} d\mu = \gen{f_3, h \cdot_c f_3} = \gen{P_c(f_3), P_c(f_3)} = 0.
\]
Thus,
\[
(\forall^\mu h) \, (\forall^\mu g) 
\; 
\gen{e_g, e_{gh}} 
= 
(\int_G F_2^{(h)} d\mu) (\int_G F_3^{(h)} d\mu) 
= 
(\int_G F_2^{(h)} d\mu) \cdot 0 
= 
0.
\qedhere
\]
\end{case}
\end{proof}

\section{A quantitative version}\label{sec:quantitative_version}

We now work out a quantitative version of the double recurrence theorem, where we consider probability groups that may not be purely mixing, but are mixing with some error (called $\e$-mixing below). 

\begin{credits}
The argument below is the same as above for the infinitary version (replacing the a.e. statements with averages), except for the proof of the approximate van der Corput lemma (\cref{approx-van-der-Corput}). The proof of the infinitary/qualitative counterpart (\cref{random_van_der_Corput}) uses a Ramsey-theoretic argument, which would still yield a quantitative bound on the error, but it would be quite rough and messy to compute. Thus, in the original version of the current paper, quantitative double recurrence was only mentioned in a remark with its proof omitted because the bound it gave was superseded by that in \cite{Austin:equidist_quasirandom}*{Theorem 1}, where a nice bound of $4 D^{-1/8}$ was obtained for $D$-quasirandom groups. However, after receiving the original version of the current paper (private communication), Austin pointed out an argument replacing the Ramsey-theoretic part of the proof with applications of the Fubini property and Cauchy--Schwarz. With Austin's permission, we use this argument to prove \cref{approx-van-der-Corput} below and obtain a slightly better bound of $3 D^{-1/4}$ for the double recurrence theorem.
\end{credits}

\subsection{Approximate mixing}

The exposition below is mainly self-contained and, although written for probability groups, the main application we have in mind is to the following class of groups:

\begin{defn}[Gowers \cite{Gowers:quasirandom_gps}]\label{def:quasirandom-groups}
For $D \ge 1$, a compact Hausdorff group $G$ is called \emph{$D$-quasirandom} if it does not admit any nontrivial unitary representations of dimension less than $D$.
\end{defn}

Below, we treat compact Hausdorff groups as probability groups as described in \examplesref{examples:prob_gps}{examples:prob_gps:compact}.

\begin{defn}[Approximate mixing]
For $\e > 0$, call a measure-preserving action $a : G \actson X$ of a probability group $(G, \mu)$ on a probability space $(X, \nu)$ \emph{$\e$-mixing} if for any $f_1, f_2 \in L^2(X,\nu)$,
\[
\int_G \big| \gen{f_1, g \cdot_a f_2}_X - \gen{P_a(f_1), P_a(f_2)}_X \big| d\mu(g)
\le 
\e \|f_1\|_{L^2} \|f_2\|_{L^2},
\]
Furthermore, call a probability group $G$ \emph{$\e$-mixing} if all of its measure-preserving actions on probability spaces are $\e$-mixing.
\end{defn}

\cite{Bergelson-Tao:mult_rec_quasirandom}*{Proposition 3}, as written, states that the right translation action of a $D$-quasirandom group on itself is $D^{-1/2}$-mixing, but running its proof for any other measure-preserving action actually yields 

\begin{prop}[Bergelson--Tao]\label{quasirandom_is_mixing}
	For all $D \ge 1$, every $D$-quasirandom compact Hausdorff group (as a probability group) is $D^{-1/2}$-mixing.
\end{prop}

\subsection{Approximate van der Corput lemma}

\begin{lemma}[Approximate van der Corput]\label{approx-van-der-Corput}
Let $\big(G, (\BC^{(k)}, \mu^{(k)})_{k \ge 1}\big)$ be a probability group and $(X,\nu)$ be a probability space. Let $(e_g)_{g \in G} \in L^2(X,\nu)^G$ be a bounded (in the $L^2$-norm) sequence such that 
\begin{enumref}{approx-van-der-Corput}
\item the function $(g,h) \mapsto \gen{e_g, e_h} : G^2 \to \C$ is $\BC^{(2)}$-measurable,

\item for every $f \in L^2(X,\nu)$, the function $g \mapsto \gen{f,e_g} : G \to \C$ is $\BC$-measurable.
\end{enumref}
For every $\e \ge 0$, if 
\[
\int_G \int_G |\gen{e_g, e_{gh}}| d\mu(g) d\mu(h) \le \e,
\]
then for all $f \in L^2(X,\nu)$,
\[
\int_G |\gen{f, e_g}| d\mu(g) \le \sqrt{\e} \ \Lnorm{f}
\]
\end{lemma}

\begin{proof}[Proof\emph{(Austin)}]
Let $\phi : G \to \C$ be defined so that $|\gen{f, e_g}| = \phi(g) \gen{f, e_g}$. Then
\begin{align*}
\int_G |\gen{f, e_g}| d\mu(g) 
&= 
\int_G \int_X \phi(g) \, f(x) \, e_g(x) \, d\nu(x) \, d\mu(g) 
\\
\eqcomment{Fubini} 
&= 
\int_X f(x) \Big(\int_G \phi(g) \, e_g(x) \, d\mu(g)\Big) \, d\nu(x) 
\\
\eqcomment{Cauchy--Schwarz} 
&\le 
\Lnorm{f} \cdot \normlr{\int_G \phi(g) e_g(\cdot) d\mu(g)}_{L^2}.
\end{align*}
But the following calculation shows that the second factor in the last term is bounded by $\sqrt{\e}$:
\begin{align*}
\normlr{\int_G \phi(g) e_g(\cdot) d\mu(g)}_{L^2}^2 
&= 
\gen{\int_G \phi(g) e_g(\cdot) d\mu(g), \int_G \phi(h) e_h(\cdot) d\mu(h)}
\\
&=
\int_X \int_G \int_G \phi(g) \, \cl{\phi(h)} e_g(x) \, \cl{e_h(x)} \, d\mu(h) \, d\mu(g) \, d\nu(x) 
\\
\eqcomment{change of variable $h \mapsto gh$} 
&= 
\int_X \int_G \int_G \phi(g) \, \cl{\phi(gh)} \, e_g(x) \, \cl{e_{gh}(x)} \, d\mu(h) \, d\mu(g) \, d\nu(x) 
\\
\eqcomment{Fubini} 
&= 
\int_G \int_G \phi(g) \, \cl{\phi(gh)} \, \gen{e_g, e_{gh}} \, d\mu(g) \, d\mu(h) 
\\
\eqcomment{triangle inequality} 
&\le \int_G \int_G |\gen{e_g, e_{gh}}| \, d\mu(g) \, d\mu(h) 
\le 
\e.
\qedhere
\end{align*}
\end{proof}

\subsection{Approximate mixing implies approximate double recurrence}

\begin{defn}[Approximate double recurrence]\label{defn:approx_dbl_rec}
	For $\e \ge 0$, call a probability group $(G,\mu)$ \emph{$\e$-doubly recurrent} if for any $f_1, f_2, f_3 \in L^2(G,\mu)$ with $L^\w$-norm at most $1$,
	\[
	\int_G 
	\left|
	\int_G f_1(x) (g \cdot_\ell f_2)(x) \, (g \cdot_c f_3)(x) \, d\mu(x) 
	- 
	\int_G f_1(x) \, P_\ell(f_2)(x) \, P_c(f_3)(x) \, d\mu(x)
	\right| 
	d\mu(g) 
	\le 
	\e.
	\]
\end{defn}

\begin{theorem}
For any $0 \le \e \le 1$, every $\e$-mixing probability group is $3 \sqrt{\e}$-doubly recurrent.
\end{theorem}
\begin{proof}
	Let $(G,\mu)$ and $f_1,f_2,f_3$ be as in \cref{defn:approx_dbl_rec} and consider the orthogonal decomposition $f = P_c(f_3) + \big(f_3 - P_c(f_3)\big)$. On one hand, \cref{mean_ergodic_theorem} implies $\norm{P_c(f_3)}_{L^\w} \le \norm{f_3}_{L^\w} \le 1$, so $\|f_3 - P_c(f_3)\|_{L^\w} \le 2$. On the other hand, Pythagorean theorem gives $\|f_3 - P_c(f_3)\|_{L^2} \le \|f_3\|_{L^2} \le \|f_3\|_{L^\w} \le 1$. Thus, noting that $e + 3 \sqrt{\e} < 3 \sqrt{e}$, our task splits into the following two:

\begin{case*}{1}
	Assuming $P_c(f_3) = f_3$ and $\norm{f_3}_{L^\w} \le 1$, prove
	\[
	\int_G 
	\left|
	\int_G f_1(x) \, f_3(x) \, (g \cdot_\ell f_2)(x) \, d\mu(x) 
	- 
	\int_G f_1(x) \, f_3(x) \, P_\ell(f_2)(x) \, d\mu(x)
	\right| 
	d\mu(g)
	\le 
	\e.
	\]
\end{case*}

\begin{case*}{2}
	Assuming $P_c(f_3) = 0$, $\norm{f_3}_{L^\w} \le 2$, and $\norm{f_3}_{L^2} \le 1$, prove
	\begin{equation}\label{eq:approx_dbl_rec_for_orthogonal}
		\int_G
		\left|
		\int_G f_1(x) \, (g \cdot_\ell f_2)(x) \, (g \cdot_c f_3)(x) \, d\mu(x)
		\right| 
		d\mu(g) 
		\le 
		\sqrt{3\e}.
	\end{equation}
\end{case*}

\noindent Case 1 is just the statement of $\e$-mixing of the left translation action applied to functions $f_1 f_3$ and $f_2$, so we focus on Case 2 now. To this end, we suppose $P_c(f_3) = 0$ and put $e_g = (g \cdot_\ell f_2)(g \cdot_c f_3)$. The approximate van der Corput lemma (\cref{approx-van-der-Corput}) reduces proving \labelcref{eq:approx_dbl_rec_for_orthogonal} to proving the following:
\[
\int_G \int_G |\gen{e_g, e_{gh}}| d\mu(g) d\mu(h) \le 3\e,
\]
For fixed $g,h \in G$, the computation done in the proof of \cref{dbl_rec} (algebraic manipulations followed by a change of variable) gives:
\[
\gen{e_g, e_{gh}} = \gen{F_2^{(h)}, g \cdot_r F_3^{(h)}},
\]
where $F_2^{(h)} \defeq f_2 (h \cdot_\ell f_2)$ and $F_3^{(h)} \defeq f_3 (h \cdot_c f_3)$. Integrating over $g$ gives:
\begin{align*}
\int_G |\gen{e_g, e_{gh}}| d\mu(g)
=& 
\int_G \left|\gen{F_2^{(h)}, g \cdot_r F_3^{(h)}}\right| d\mu(g) 
\\
\eqcomment{triangle inequality} 
\le& 
\int_G \left|\gen{F_2^{(h)}, g \cdot_r F_3^{(h)}} - \gen{P_r(F_2^{(h)}), P_r(F_3^{(h)})}\right|d\mu(g)
\\ 
&+ 
\left|\gen{P_r(F_2^{(h)}), P_r(F_3^{(h)})}\right| 
\\
\eqcomment{right translation is $\e$-mixing} 
\le& 
\, \e \Lnorm{F_2^{(h)}} \Lnorm{F_3^{(h)}}
+
\left|\gen{P_r(F_2^{(h)}), P_r(F_3^{(h)})}\right|.
\end{align*}

\noindent 
But $\Lnorm{F_2^{(h)}} \le \norm{F_2^{(h)}}_{L^\w} \le 1$ and 
$
\Lnorm{F_3^{(h)}} \le \norm{f_3}_{L^\w} \Lnorm{h \cdot_\ell f_3} = \norm{f_3}_{L^\w} \Lnorm{f_3} \le 2.
$
As for the last term, because right multiplication is transitive, $P_r(f) \equiv \int_G f$ $\mu$-a.e. for any $f \in L^2(G,\mu)$, so
\begin{align*}
\left|\gen{P_r(F_2^{(h)}), P_r(F_3^{(h)})}\right|
&= 
\left|\int_G F_2^{(h)}\right| \left|\int_G F_3^{(h)}\right| 
\\
&\le  
\norm{F_2^{(h)}}_{L^\w}^2 \cdot |\gen{f_3, (h \cdot_c f_3)}|
\\
&\le 
|\gen{f_3, (h \cdot_c f_3)}|.
\end{align*}
Finally, putting it all together and integrating over $h$ gives:
\begin{align*}
	\int_G \int_G |\gen{e_g, e_{gh}}| \, d\mu(g) \, d\mu(h)
	&\le
	2\e + \int_G \left| \gen{f_3, (h \cdot_c f_3)} \right| d\mu(h)
	\\
	\eqcomment{$P_c(f_3) = 0$}
	&=
	2\e + \int_G \left| \gen{f_3, (h \cdot_c f_3)} - \gen{P_c(f_3),P_c(f_3)} \right| d\mu(h)
	\\
	\eqcomment{conjugation is $\e$-mixing}
	&\le
	2\e + \e
	=
	3\e.
	\qedhere
\end{align*}
\end{proof}

\cref{quasirandom_is_mixing} and the last theorem give \cite{Austin:equidist_quasirandom}*{Theorem 1} with a slightly better bound:

\begin{cor}
For all $D \ge 1$, every $D$-quasirandom compact Hausdorff group (as a probability group) is $3 D^{-1/4}$-doubly recurrent.
\end{cor}

\begin{acknowledgements}
This paper owes a great deal to \Slawek Solecki: his short note explaining the main result in \cite{Bergelson-Tao:mult_rec_quasirandom}, as well as generalizing Proposition 3 of \cite{Bergelson-Tao:mult_rec_quasirandom} to arbitrary actions of quasirandom groups on finite sets, was what cleared things up for me and gave the right prospective. Also, many thanks to Evgeny Gordon for useful conversations and comments; in particular, for bringing to my attention the fact that a notion of probability groups had already been defined in \cite{Weil}, as well as for pointing out that the Loeb measure is \emph{countably additive} even when defined on an ultraproduct of \emph{finitely additive} probability spaces. Further thanks to Terence Tao for very useful comments and suggestions, as well as to Julien Melleray and Benjamin Weiss for pointing out some errors in my original definitions of probability groups and their actions. Also, thanks to Asgar Jamneshan for pointing out the measurability of multiplication issue in non-second-countable compact groups. Many thanks to the anonymous referee for carefully reading the paper and suggesting the example of the ultraproduct of amenable unimodular locally compact groups. Finally, I am grateful to Tim Austin for allowing me to use his argument in the approximate van der Corput lemma, as well as for his other comments and corrections.
\end{acknowledgements}


\bigskip



\end{document}